\documentclass[11pt, a4paper]{article}  
\usepackage{helvet}	             
\usepackage{geometry}
\usepackage{multirow}
\usepackage{amsmath}
\usepackage{amssymb}
\usepackage[usenames]{color}
\usepackage{graphicx}
\usepackage{amsthm}
\usepackage[dvips]{epsfig}
\usepackage{booktabs}

\newcommand\C{{C}}	
\newtheorem{theorem}{Theorem}
\newtheorem{definition}{Definition}
\newtheorem{corollary}{Corollary}
\newtheorem{remark}{Remark}
\newtheorem{problem}{Problem}

\begin{document}

\begin{center}
  {\Large \bf
A family of extremum seeking laws
\vspace{.25em}
 for a unicycle model  with a moving target:
 \vspace{.25em}
  theoretical and experimental studies$^{*}$}
\end{center}
\vspace{0.5em}
\begin{center}
  { \bf
  Victoria Grushkovskaya$^{1,3}$, Simon Michalowsky$^{1}$, Alexander Zuyev$^{2,3}$,\\ {Max May}$^{1}$ and Christian Ebenbauer$^{1}$
}
\end{center}
\vspace{0.5em}
\begin{center}
\scriptsize  $^{1}$Institute for Systems Theory and Automatic Control, University  of Stuttgart, 70569 Stuttgart, Germany
        \{grushkovskaya,michalowsky,may,ce\}@ist.uni-stuttgart.de\\
        $^{2}$Max Planck Institute for Dynamics of Complex Technical Systems, 39106 Magdeburg, Germany\\ zuyev@mpi-magdeburg.mpg.de\\
        $^{3}$Institute of Applied Mathematics and Mechanics, National Academy of Sciences of Ukraine,  \\ 84100 Sloviansk, Ukraine
\end{center}
\let\thefootnote\relax\footnote{$^{*}$ This paper is an extended version of~\cite{GMZME18} which contains the proof of the main result. 
This work was supported in part by the German Research Foundation (DFG, EB 425/4-1).
}
\begin{abstract}
In this paper, we propose and practically evaluate a class of gradient-free  control functions ensuring the motion of a unicycle-type system towards the extremum point of a time-varying cost function.
We prove that the unicycle is able to track the extremum point, and illustrate our results by numerical simulations and experiments that show that the proposed control functions exhibit an improved tracking performance in comparison to standard extremum seeking laws based on Lie bracket approximations.
\end{abstract}

\section{Introduction}

{Extremum seeking} typically refers to the problem of constructing a gradient-free control law that ensures the motion of a dynamical system to the minimum (or maximum) of a partially or completely unknown and possibly time-varying cost or performance function.
Over the past decades, significant advances in the theory and applications of extremum seeking have been made, see, e.g., \cite{Krst03,review}.
Today, there exist many ways to design and analyze extremum seeking laws exploiting, e.g., averaging and singular perturbation techniques, Lie bracket approximation techniques, least squares estimation approaches, stochastic and hybrid approaches, see, e.g., \cite{Att-Joh-Gus-15,Due-Krs-Sch-Ebe-15,DurrAuto,Guay03,Har-Wou-Nes-13,Khong15,Kr00,Liu12,Nes10,Pov-Tee-17}.
Many extremum seeking schemes use control functions
{depending on the current value of the cost function modulated by}
time-periodic oscillating excitation (or dither, learning) signals in order to explore and extract sufficient information from the dynamical system and/or from the unknown cost function to solve the extremum seeking problem.

The choice of the {control function} 
as well as the excitation signals plays an important role for the performance of the extremum seeking scheme~\cite{Chi07,Nes09,Sch16,TNM08}.
In the recent paper~\cite{GZE}, a broad family of control functions for extremum seeking schemes based
on Lie bracket approximations was presented for systems with single-integrator dynamics and for time-invariant cost functions. This class of controls has several favorable
properties including the possibility of adapting and constraining the amplitude of the excitation
signal.
Moreover, in the paper~\cite{GDEZ17} the extremum seeking problem for time-varying cost functions
has been considered in the framework of Lie bracket approximations, but again with single-integrator dynamics and {with standard  control functions as used in \cite{DurrAuto}}.
{The first {contribution} of this paper is a whole family of control functions which enables a system with unicycle-type dynamics to approximate the gradient-like flow of a time-varying cost functions. This result justify the use of gradient-free controllers presented in \cite{GZE}  in time-varying extremum seeking problems.
  For the sake of simplicity, we consider a distance-like time-varying cost function. Such problems arise, for example, when a robot has to follow a moving target {(tracking problem)} and only the distance (but not the relative position) to the moving target can be measured.
Although gradient-free control laws for the tracking of a moving target  have been previously considered  (see, e.g., ~\cite{De05,Due-Sta-Joh-Ebe-12,Hua13,Mand15,Mos15,Sah12,Schei12,Sch14,Vweza2015,Yu05,Zhu13}), the main advantage of the proposed family of extremum seeking laws is, on the one hand,
the high flexibility in designing the {control functions} such that they meet
further specifications like input constraints, and,
on the other hand, the family of control functions ensures rigorous stability and tracking properties.

As it will be shown, some important control strategies in the proposed class are not continuously differentiable. In view of this, the second {contribution} of this paper is {the relaxation of} the ``$\C^2$-requirement'' for the Lie bracket approximation approach~\cite{DurrAuto,GDEZ17}. Instead, we will require the continuity of Lie derivatives. This result will allow {us} to exploit a much  wider class of admissible extremum seeking laws. Extremum seeking systems with non-$\C^2$ vector fields were also considered in~\cite{GZE,Sch14b} for time-invariant cost functions. However,  the results of the above papers are not directly applicable for time-varying extremum seeking problems.}

{{As the third} contribution} of this paper, {we}  show by numerical simulations and experiments
with {a mobile robot} that the high flexibility of
 the proposed control functions can be utilized to significantly improve the tracking behavior { in comparison to standard extremum seeking approaches considered, for example, in~\cite{DurrAuto}.}

The rest of this paper is organized as follows. In Section~\ref{Prel}, we formulate the extremum seeking problem and recall some results on the Lie bracket approximation approach.
 Section~\ref{Main} presents a class of extremum seeking laws for a unicycle-type system and stability results  for time-varying cost functions.
The numerical simulations and experiments for several extremum seeking laws with different qualitative properties are discussed in Section~\ref{appl}.

\section{Problem Statement and Preliminaries}~\label{Prel}
\vspace{-2em}
\subsection{Problem statement}
Consider a unicycle model
\begin{equation}\label{uni}
\begin{aligned}
&\dot x_1=u\cos(\Omega t),\\
&\dot x_2=u\sin(\Omega t),\\
\end{aligned}
\end{equation}
where {{$x=(x_1,x_2)^T\in\mathbb R^2$ is the state, $u\in\mathbb R$ is the control}}, and $\Omega>0$ is constant angular velocity.
These equations correspond  to the standard unicycle model where the equation for the angular velocity
$\dot \theta = \Omega$, $\theta(0)=0$, has been eliminated (see, e.g., \cite{Sch14,GE16} for details).  Note that the model with non-constant angular velocity can also be considered
using, e.g., singular perturbations techniques~\cite{Due-Krs-Sch-Ebe-15}.
In the sequel, $\mathbb R^+$ denotes the set of all non-negative real numbers.
In this paper, we address  the following problem:

\begin{problem}
\textit{For a cost function $\hat J\in  C^2(\mathbb R^2 \times \mathbb R^2; \mathbb R)$,  $\hat J{=}\hat J({x},\gamma)$,
with an (unknown)  $\gamma:\mathbb R^+{\to}\mathbb R^2$, $\gamma=(\gamma_1(t),\gamma_2(t))^T$ and a unique (possibly time-varying) minimum
  $x^*(\gamma(t))$, the goal is to construct a control function $ u = u\big(t,\hat J(x(t),\gamma(t))\big)$  that  asymptotically steers  system~\eqref{uni} to an arbitrary small neighborhood of $x^*(\gamma(t))$ as $t\to +\infty$.
}\end{problem}
In particular, if
 \begin{equation}\label{J}
 \hat J(x,\gamma)=J(x-\gamma)= \kappa\|x-\gamma\|^2 \text{ with some }\kappa>0,
 \end{equation}
 then the above extremum seeking problem leads to steering the control system~\eqref{uni} to an arbitrary small neighborhood of {the curve } $x^*(t)=\gamma(t)$.
\subsection{{Notations}}
For  $f_i,f_j:\mathbb R\times\mathbb R^n\to\mathbb R^n $, $x\in\mathbb R^n$, we denote the Lie derivative $L_{f_i}f_j(t, x)=\lim\limits_{y\to0}\frac{1}{y}(f_j(t,x+yf_i(t,x))-f_j(t,x))$ and the Lie bracket $[f_i,f_j](t, x){=}L_{f_i}f_j(t, x){-}L_{f_j}f_i(t,x)$;
for a function $h\in {{C}^1}(\mathbb R^n\times \mathbb R^k;\mathbb R)$ and $\xi\in\mathbb R^n$, the gradient of $h$ with respect to $\xi$ is denoted as
$\nabla_\xi h(\xi,\eta)=\frac{\partial h(\xi,\eta)}{\partial \xi}^T$;
for  a one-parameter family of non-empty sets $\mathcal L_t\subset\mathbb R^n$, \linebreak$t\in\mathbb R^+$,
 $
 B_{\delta}(\mathcal L_{t}){=}\cup_{y{\in}\mathcal L_{t}}\{x{\in}\mathbb R^n:\|x{-}y\|{<}\delta\}
 $ is a $\delta$-neigh\-borhood of the set
$\mathcal L_{t}$ at time $t$, and
 $\bar B_{\delta}(\mathcal L_{\lambda,t})$ denotes the closure of $B_{\delta}(\mathcal L_{\lambda,t})$.
{For $m,n\in\mathbb Z$, the notation $i=\overline{m,n}$ means that $i$ takes all possible values from the set $\{m,m+1,\dots,n\}$.}

  \subsection{Preliminaries}
In this section, we summarize some basic facts about Lie bracket approximation techniques.
Consider a control-affine system
\begin{equation}\label{aff}
  \dot x=f_0(t,x)+\sum_{j=1}^{{\ell}} f_j(t,x)\sqrt\omega u_j(\omega t),
\end{equation}
where {$x{=}(x_1(t),\dots,x_n(t))^T{\in}D\subseteq\mathbb R^n$}, $x(t_0){=}x^0{\in}D$, $t_0{\in}{\mathbb R}^+$, $\omega{>}0$.
Let the following assumptions be satisfied: \\
 (A1)~$f_i\in C^2(\mathbb R^+\times D;\mathbb R^n)$, $i=\overline{0,\ell}$;\\
(A2)~the functions $\|f_i(t,x)\|$,$\|\frac{\partial f_i(t,x)}{\partial t}\|$,$\|\frac{\partial f_i(t,x)}{\partial x}\|$,$\|\frac{\partial^2 f_j(t,x)}{\partial t\partial x}
\|$, $\left\|\frac{\partial [f_j,f_k](t,x)}{\partial t}\right\|$,
      $\left\|\frac{\partial [f_j,f_k](t,x)}{\partial x}\right\|$ are bounded on each compact set $x{\in}\mathcal X{\subseteq}D$
     uniformly in $t{\ge}0$, for  {{$i{=}\overline{0,\ell}$, $j{=}\overline{1,\ell}$, $k{=}\overline{j,\ell}$.}} \\
(A3)~the functions $u_j$ are Lipschitz continuous and $T$-periodic with some $T{>}0$, and $\int_0^Tu_j(\tau)d\tau=0$, $j=\overline{1,\ell}$.\\
Since in time-varying extremum seeking problems it is often necessary to investigate the stability of a family of sets (instead of a single set), we will make use of the following definitions which can be found, e.g., in~\cite{GDEZ17}.
 \begin{definition}
A family of  {non-empty} sets   {$\mathcal L_{t} \subseteq D $, $ t \in{\mathbb R}^+$,} is said to be \textit{locally practically uniformly asymptotically stable} for~\eqref{aff} if it is\\
$\bullet$ \textit{practically uniformly stable}: for any $\varepsilon{>}0$ there exist $\delta,\omega_0{>}0$ such that, for all $t_0{\in}\mathbb R^+$ and $\omega{>}\omega_0$, the following property holds {for the solutions of~\eqref{aff}}:
        $$x^0{\in} B_{\delta}(\mathcal L_{t_0})\Rightarrow x(t){\in}  B_{\varepsilon}(\mathcal L_{t})\text{ for all }t{\ge} t_0;$$
$\bullet$  $\hat\delta$-\textit{practically uniformly attractive}  with some $\hat\delta>0$: for every $\varepsilon>0$ there exist $t_1\ge0$ and $\omega_0>0$ such that, for all $t_0\in\mathbb R^+$ and $\omega>\omega_0$, the following property holds {for the solutions of~\eqref{aff}}:
        $$x^0{\in} B_{\hat\delta}(\mathcal L_{t_0})\Rightarrow x(t)\in B_{\varepsilon}(\mathcal L_{t})\text{ for all }t{\ge} t_0+t_1;$$
$\bullet$ the solutions of~\eqref{aff} are \textit{practically uniformly bounded}: for each $\delta{>}0$ there are $\varepsilon{>}0$ and $\omega_0{>}0$ such that, for
    all $t_0{\in}\mathbb R^+$ and $\omega{>}\omega_0$, the following property holds {for the solutions of~\eqref{aff}}:
          $$x^0{\in} B_{\delta}(\mathcal L_{t_0})\Rightarrow x(t){\in}  B_{\varepsilon}(\mathcal L_{t})\text{ for all }t{\ge} t_0.$$
  If the attractivity property holds for every $\hat\delta{>}0$, then the family of sets $\mathcal L_{t}$ is called  \textit{semi-globally practically uniformly asymptotically stable} for~\eqref{aff}.
 For systems independent {of} $\omega$ we omit the terms \textit{practically} and \textit{semi}.
 \end{definition}
  The following result from~\cite{GDEZ17} allows to establish practical asymptotic stability properties of~\eqref{aff} from asymptotic stability properties of the so-called Lie bracket system.
{\begin{theorem}~\label{ourthm}
  \textit{ Let {{(A1)--(A3)}} hold. Suppose that a family of  sets $ \mathcal L_t\subset D$ is locally (globally) uniformly asymptotically stable for the Lie bracket system
  \begin{equation}\label{affLie}
    \dot {\bar x}{=}f_0(t,\bar x){+}\frac{1}{T} {\sum_{i=1}^n \sum_{\substack{ j=2 \\ j > i}}^{n} } [f_i,f_j](t,\bar x){\int_0^{T}}{{\int_0^\theta}{u_j(\theta)u_i(\tau)}d\tau}d\theta,
  \end{equation}
 and suppose that there exists a compact set  $S{\subset}\mathbb R^n$ such that  $\mathcal L_t{\subseteq }S$ for all $t\in\mathbb R^+$.
   Then $ \mathcal L_t$ is locally (semi-globally) practically uniformly asymptotically stable for~\eqref{aff}.}
  \end{theorem}	}
  In order to characterize the stability and tracking behavior of~\eqref{uni} with respect to
time-varying cost functions {of the form}~\eqref{J}, we will consider {the following} family of level sets of the cost function $\hat J$:
$$
\mathcal L_{\lambda,t}=\{x\in\mathbb R^2: \hat J(x,\gamma(t)) - \hat J(x^*(\gamma(t)),\gamma(t))\le \lambda\},{\; \lambda,t\geq0}.
$$

Let us emphasize that  most of {the} results presented in the following {apply} to general
time-varying cost functions not necessarily restricted to the form~\eqref{J}. {However, in
target tracking applications, the cost function often takes the form~\eqref{J} with a possibly unknown~$\gamma$.	}
\section{Main results}~\label{Main}
\vspace{-2em}
\subsection{Family of extremum seeking controls}~\label{sec_uni}
{In~\cite{GZE}, we have introduced a novel family of extremum seeking controls for systems with integrator dynamics
and time-invariant cost {functions}. In this section, we show that a similar result can be obtained for system~\eqref{uni} with time-varying cost functions.}
\begin{theorem}~\label{classthm}
\textit{Consider the control function
\begin{equation}\label{cont}
\begin{aligned}
  u^\omega(t,\hat J(x,\gamma))=\sqrt{\vartheta\alpha\omega}\big(&F_1(\hat J(x,\gamma)) {\cos(\omega t)}+F_2(\hat J(x,\gamma)){\sin(\omega t)} \big),
\end{aligned}
\end{equation}
where $\omega=k\Omega$ with some $k\in\mathbb N$, {$k>1$}, $\alpha=4(1-k^{-2})$,
$\vartheta>0$, and the functions $F_1,F_2$ are such that $F_s\in C(\mathbb R;\mathbb R)$, $F_s\circ \hat J\in C^1(\mathbb R^n \times \mathbb R^n ;\mathbb R)$ ($s=1,2$), and
\begin{equation}\label{class}
F_2(z)=-F_1(z)\int{{\frac{dz}{F_1(z)^2} } }.
\end{equation}
 Then the Lie bracket system corresponding to the closed-loop system~\eqref{uni} with control function~\eqref{cont} has the form
  \begin{equation}\label{uni_Lie}
\begin{aligned}
{\dot{\bar x}}=-{{\vartheta}}\nabla_{\bar x} \hat J(\bar x,\gamma)+\Phi(\bar x,\gamma),
\end{aligned}
 \end{equation}
 with $\Phi(\bar x,\gamma)=(-\varphi_2(\bar x,\gamma),\varphi_1(\bar x,\gamma))^T$,}
$$
 \varphi_s(\bar x,\gamma)=\frac{1}{2k}\nabla_{\bar x_s}\Big(F_1^2(\hat J(\bar x,\gamma))+F_2^2(\hat J(\bar x,\gamma))\Big),\;s=1,2.
$$
\end{theorem}
\vspace{1em}
\begin{remark}
  In formula~\eqref{class}, we assume that $F_1(z)\ne0$ except for at most a countable set of isolated zeros $Z^*=\{z_k^*\}$.
  We treat the function $\Psi_1(z):=\int{{\frac{dz}{F_1(z)^2} } }$ as an antiderivative of $\displaystyle\frac{1}{F_1(z)^2}$ defined on the open set ${\mathbb R}\setminus Z^*$, so that~\eqref{class} holds as an identity with continuous functions in a neighborhood of each point $z\notin Z^*$.
  As the functions $F_1$ and $F_2$ are assumed to be globally continuous, formula~\eqref{class} is treated
  in the sense that $F_2(z_k^*)=-\lim\limits_{z\to z_k^*}F_1(z)\Psi_1(z)$ at each $z_k^*\in Z$.
\end{remark}
\vspace{1em}
\emph{Proof of Theorem~2:}
Substituting the controls~\eqref{cont} into system~\eqref{uni}, we obtain a system of the
{form}~\eqref{aff}
{given by}
\begin{align}~\label{uninew}
 &\quad\quad\quad\dot x=\sum_{j=1}^4 f_j(t,x)\sqrt\omega v_j(\omega t),
\end{align}
{where} {$v_j$} are the new {$ \frac{2\pi k}{\omega}$-periodic} inputs,
\begin{align*}
&v_1(\omega t)=\cos(\omega t)\cos\Big(\frac{\omega t}{k}\Big),\,
v_2(\omega t)=\sin(\omega t)\cos\Big(\frac{\omega t}{k}\Big),\\
&v_3(\omega t)=\cos(\omega t)\sin\Big(\frac{\omega t}{k}\Big),\,
v_4(\omega t)=\sin(\omega t)\sin\Big(\frac{\omega t}{k}\Big),
\end{align*}
and the new vector fields are $
 f_i(t,x){=}\big(F_i(\hat J(x,\gamma(t)),0\big)^T,$
$f_{i+2}(t,x){ =}\big(0,F_i(\hat J(x,\gamma(t))\big)^T,\ i=1,2$. Direct construction of system~\eqref{affLie} with the use of~\eqref{class} completes the proof.\quad\quad\quad
$\blacksquare$
\\
Notice that the above  control laws leave a lot of freedom for tuning by choosing
the functions $F_1$ and $F_2$ appropriately.

\subsection{{Stability conditions}}~\label{sec_gen}
{ If the functions $F_1$ and $F_2$ are of class $C^2$, then the stability properties of the unicycle model~\eqref{uni} controlled by~\eqref{cont} can be deduced from the stability properties of the corresponding Lie bracket system~\eqref{uni_Lie}. This directly follows from Theorem~\ref{ourthm}:
\
\begin{corollary}
     \textit{Let  the functions $F_s(\hat J(\cdot,\cdot))$ satisfy (A1)--(A2), and suppose that a one-parameter family of  sets $ \mathcal L_{\lambda,t}$ is locally (globally) uniformly asymptotically stable for~\eqref{uni_Lie} with some $\lambda>0$, and there exists a compact set  $S{\subset}\mathbb R^n$ such that  $\mathcal L_{\lambda,t}{\subseteq }S$ for all $t\in\mathbb R^+$,
   then $ \mathcal L_{\lambda,t}$ is locally (semi-globally) practically uniformly asymptotically stable for
	\eqref{uni} with controls~\eqref{cont}.}
\end{corollary}
In combination with asymptotic stability conditions of families of sets for the Lie bracket system, the above result describes a solution to Problem~1 with a wide class of time-varying cost functions $\hat J$, provided that $\lambda$ is small enough. Although, unlike~\cite{GDEZ17}, the Lie bracket system for~\eqref{uni} is not the exact gradient flow of the cost function, it admits the same asymptotic stability conditions for $\mathcal L_{\lambda,t}$ as proposed in~\cite{GDEZ17} because of the property $(\Phi(\bar x,\gamma),\nabla \hat J(\bar x,\gamma))\equiv 0$.
However, many functions $F_1,F_2$ described by~\eqref{class}   fail to satisfy the $\C^2$-condition at the origin (so that  Corollary~1 and the results of~\cite{GDEZ17} are not applicable), but exhibit much better performance in comparison with systems with smooth vector fields  (see~\cite{Sch14b,GZE,SD17} and Section~~IV for some examples). To overcome such limitation,  we will present stability results under  relaxed assumptions.
 Note that, although extremum seeking problems for systems with non-$C^2$ vector fields were previously considered, e.g., in~\cite{Sch14b,GZE}, the results of the above papers are not applicable because of several reasons. First, it is easy to see that the time-varying function {$x^*(\gamma(t))$} is not a solution of system~\eqref{uni_Lie}, {therefore, the considered problem cannot be reduced to control design in a neighborhood of an admissible trajectory}.
 Second, the approximation result  in~\cite{Sch14b} has been proved under the assumption that the Lie bracket system possesses $C^2$-vector fields: $[f_i,f_j]\in C^2(\mathbb R^+\times D)$. However, the function $\varphi_1,\varphi_2$ in Theorem~\ref{classthm} do not necessary satisfy this requirement. The following result establishes the stability of the unicycle model~\eqref{uni} controlled by~\eqref{cont} under relaxed assumptions.
For  clarity of presentation and because of  space limitations,
our next theorem and its proof will be stated for  system~\eqref{uni} with time-varying cost functions of the form~\eqref{J}.
 It is expected that similar
results can be obtained for a wide class of time-varying cost functions (but with more
involved conditions). We leave the general case for future studies.}
\begin{theorem}~\label{thm_quad}
\textit{Let the cost $J=J(x-\gamma(t))$ be of the form~\eqref{J},
$\rho>0$,
$D= \displaystyle\cup_{t\ge 0}\mathcal L_{\rho,t}$, and $F_1,F_2$ be defined from~\eqref{class}. \\Assume that:  \\
$(B1)$~$F_i\circ J \in C^2(D\setminus\{0\};\mathbb R)$,
$L_{F_j}(F_i\circ J)\in C(D;\mathbb R)$, and $L_{F_k}L_{F_j}(F_i\circ J)\in C(D;\mathbb R)$ for all $i,j,k\in\{1,2\}$;\\
$(B2)$~the first-order partial derivatives of $F_i\circ J$ and of $L_{F_i}(F_j\circ J)$ are uniformly bounded in $D\setminus\{0\}$ for all $i,j,k\in\{1,2\}$;\\
$(B3)$~$\gamma\in C^1(\mathbb R^+;\mathbb R^2)$, and  there exists a $\nu>0$ such that $\|\dot\gamma(t)\|\le\nu$  for all $t\in\mathbb R^+$.\\
Then, for any $\lambda\in(0,\rho)$, $\delta\in\Big(0,\frac{\sqrt\rho-\sqrt\lambda}{\sqrt\kappa}\Big)$, and $\vartheta>\frac{\nu}{2\sqrt{\kappa\lambda}}$,
the family of sets
\begin{equation}\label{set}
\mathcal L_{\lambda,t}=\{x\in\mathbb R^2: J(x-\gamma(t)) \le \lambda\},\;t\in\mathbb R^+
\end{equation}
 is  practically uniformly {asymptotically} stable for system~\eqref{uni} with  $x^0\in B_{\delta}(\mathcal L_{\lambda,t_0})$.
}
 \end{theorem}

The proof is in Appendix~A. Note that the proof technique is similar to~\cite[Theorem 3]{GZE}. However, since the results of~\cite{GZE} are proved for the case of {time-invariant} vector fields $f_i(x)$ and constant $x^*$, they are not directly applicable. The proof of Theorem~\ref{thm_quad} requires some extensions of the approach of~\cite{GZE} to control-affine systems with time-varying vector fields and non-vanishing drift term. Furthermore, unlike many other results on the time-varying extremum seeking problems, we do not assume  that $\gamma(t)$ is uniformly bounded.

\section{Numerical Simulations and Experiments}~\label{appl}
In this section, we illustrate our results with examples and discuss some interesting choices of the functions $F_1$ and $F_2$ in the control law~\eqref{class}.
\subsection{Moving target tracking}
\begin{figure*}[t]
     {\centering
      \includegraphics[width=0.32\linewidth]{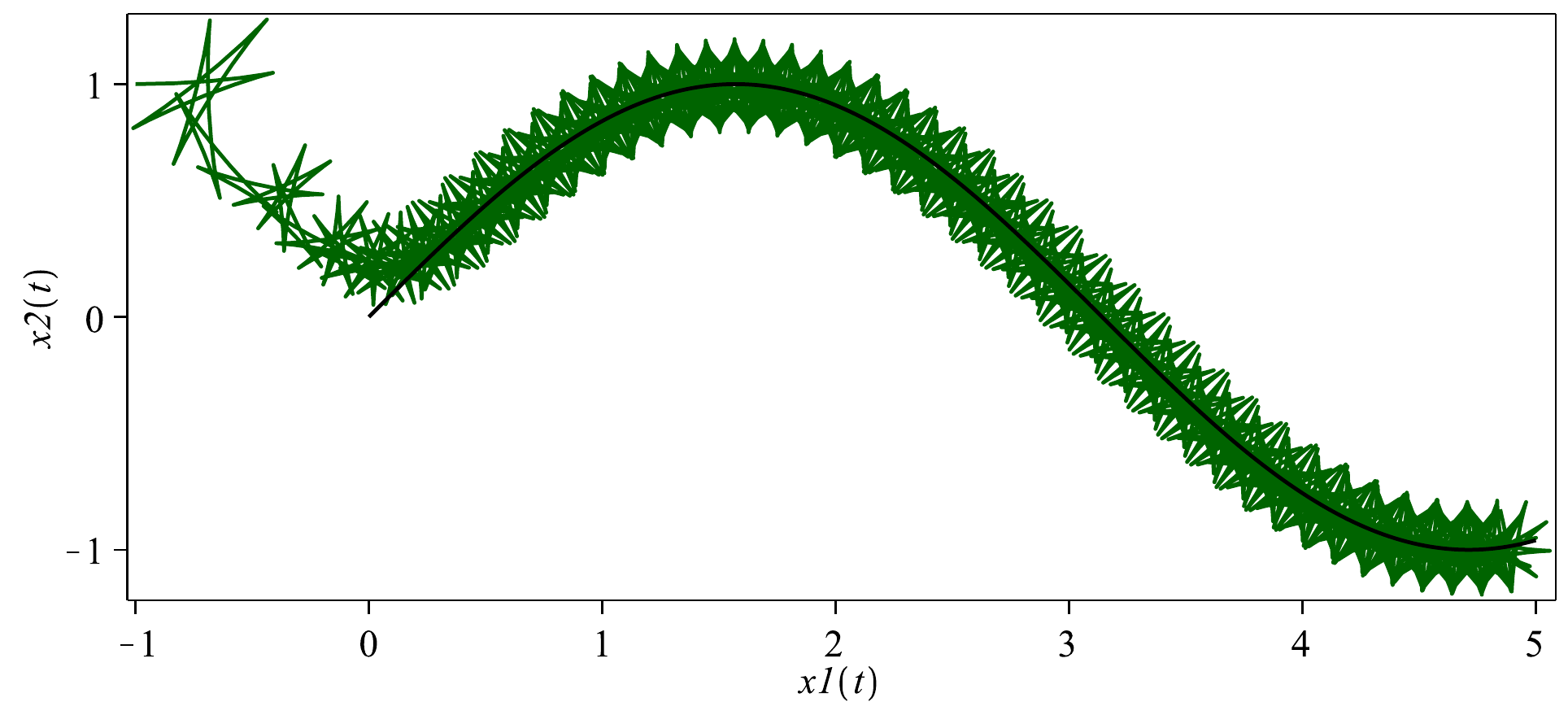} \includegraphics[width=0.32\linewidth]{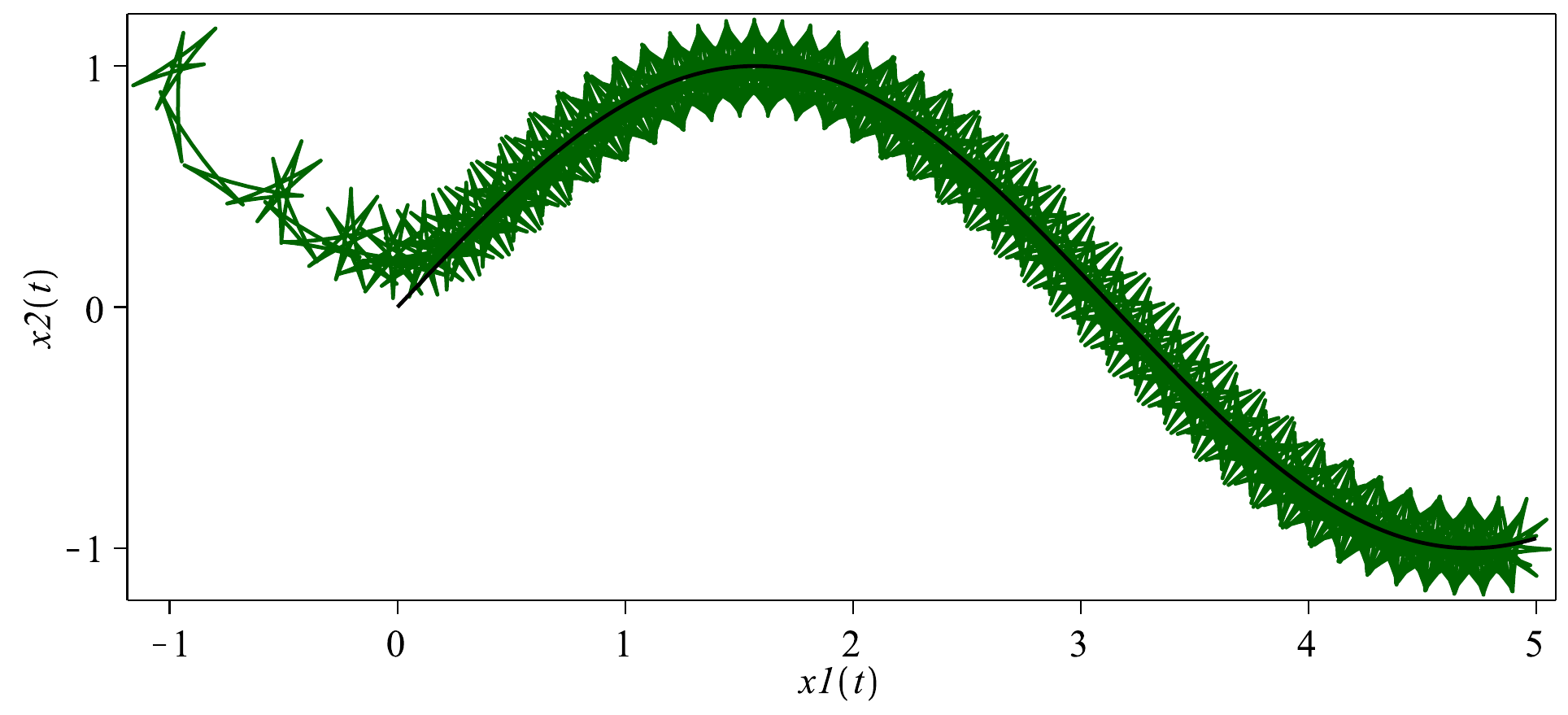} \includegraphics[width=0.32\linewidth]{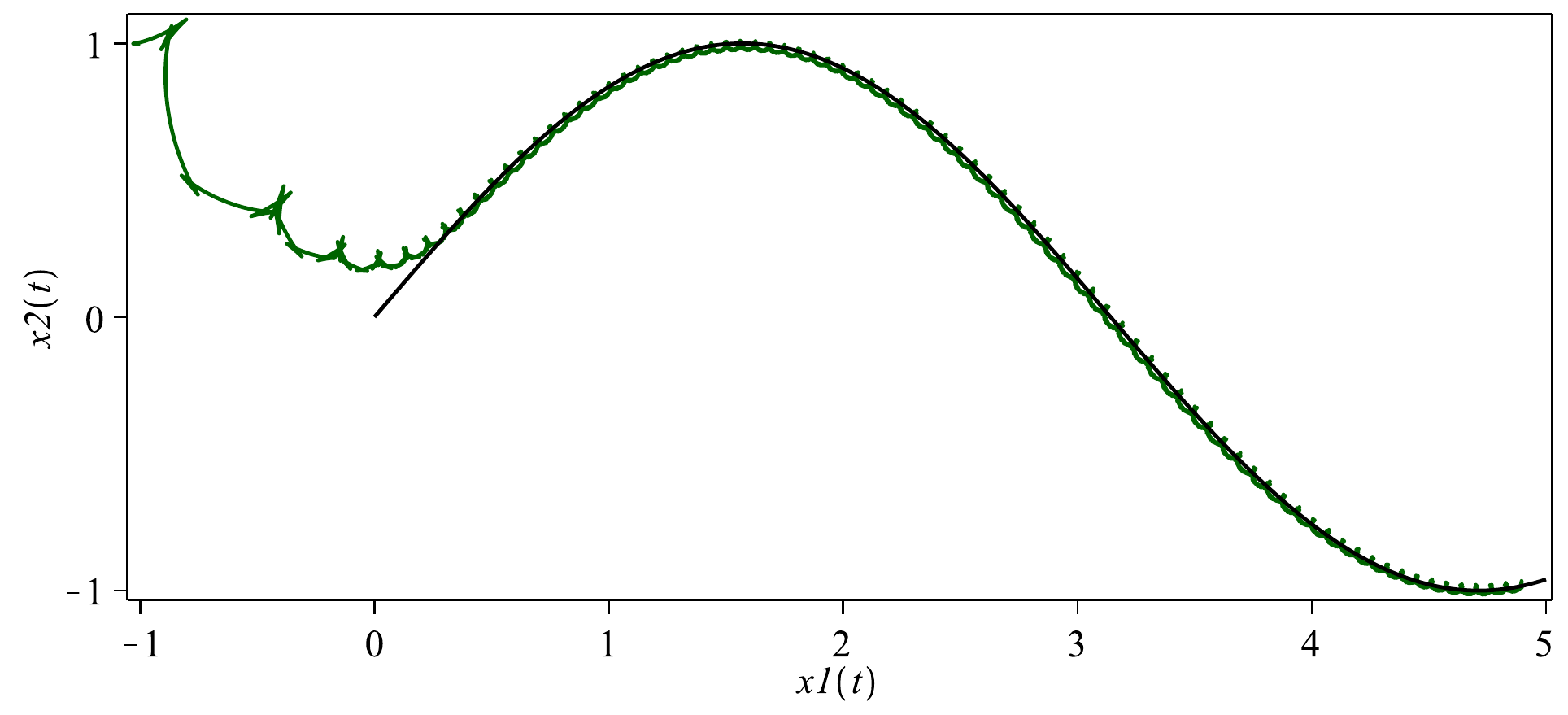}\\
      \includegraphics[width=0.32\linewidth]{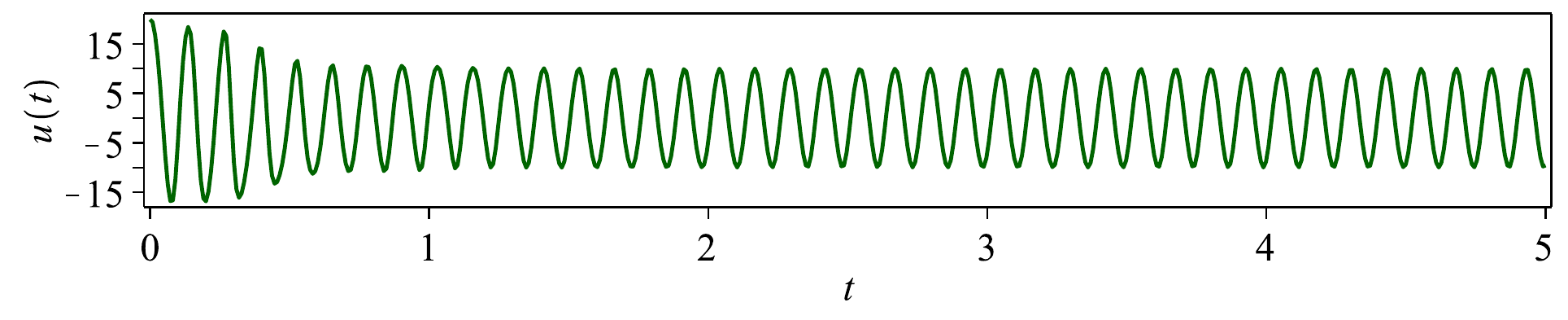}  \includegraphics[width=0.32\linewidth]{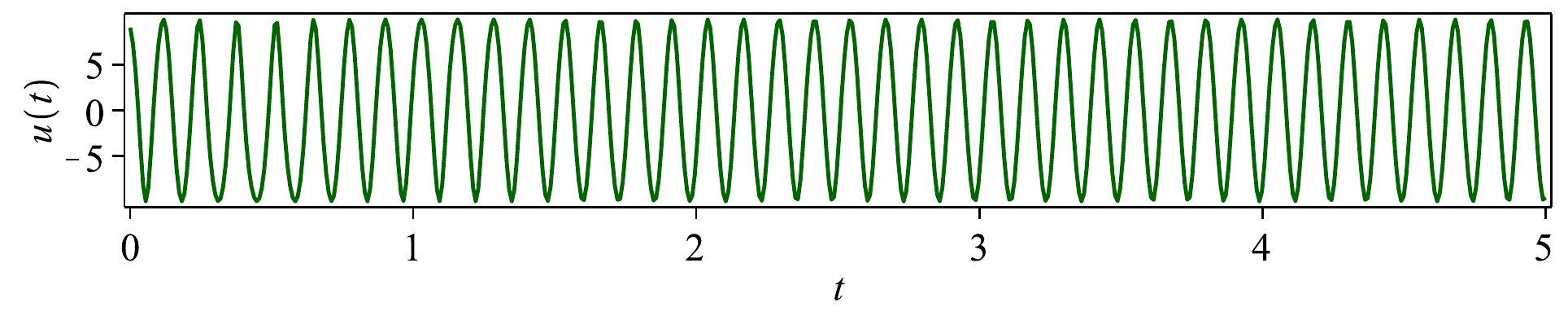} \includegraphics[width=0.32\linewidth]{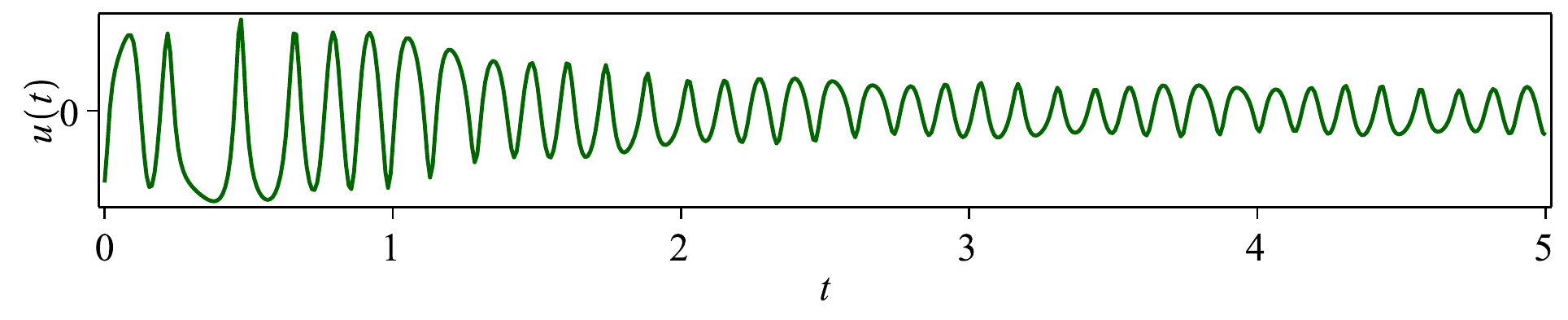}\\
      \includegraphics[width=0.32\linewidth]{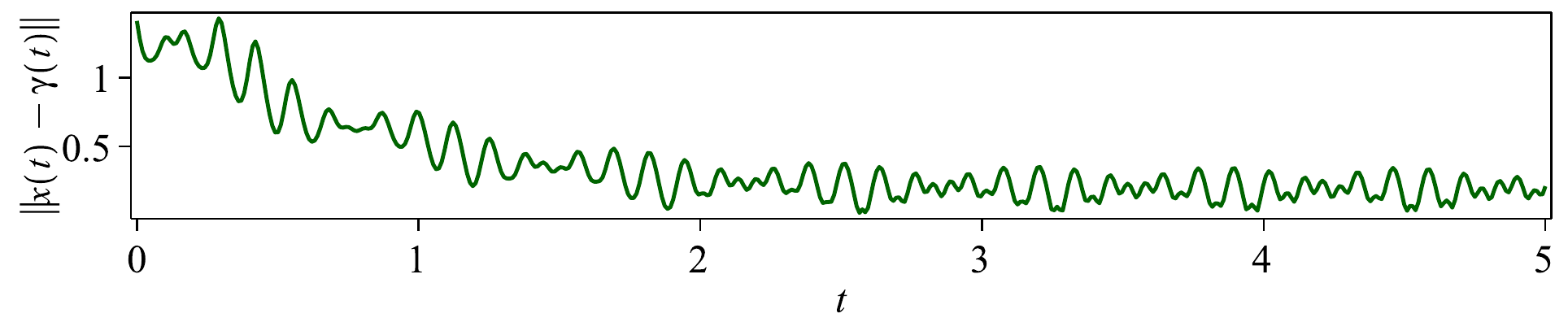}  \includegraphics[width=0.32\linewidth]{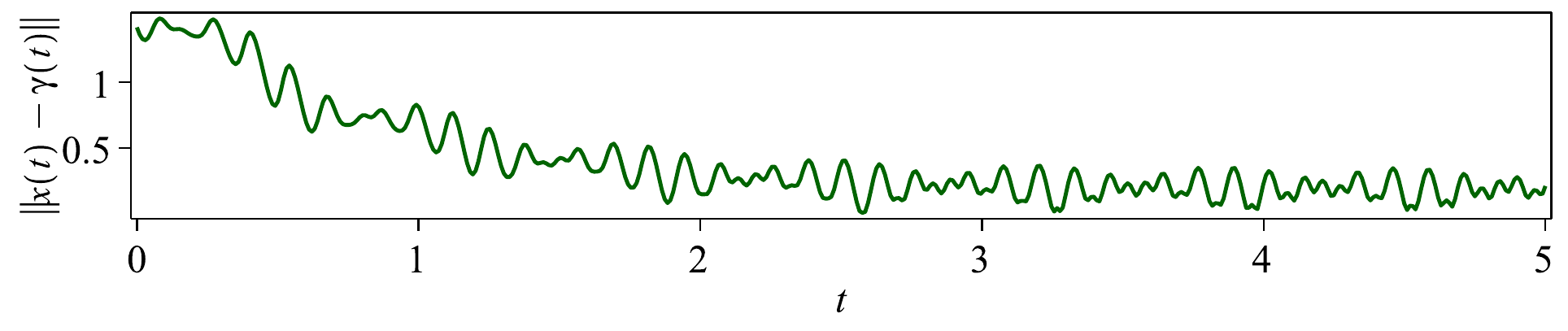}  \includegraphics[width=0.32\linewidth]{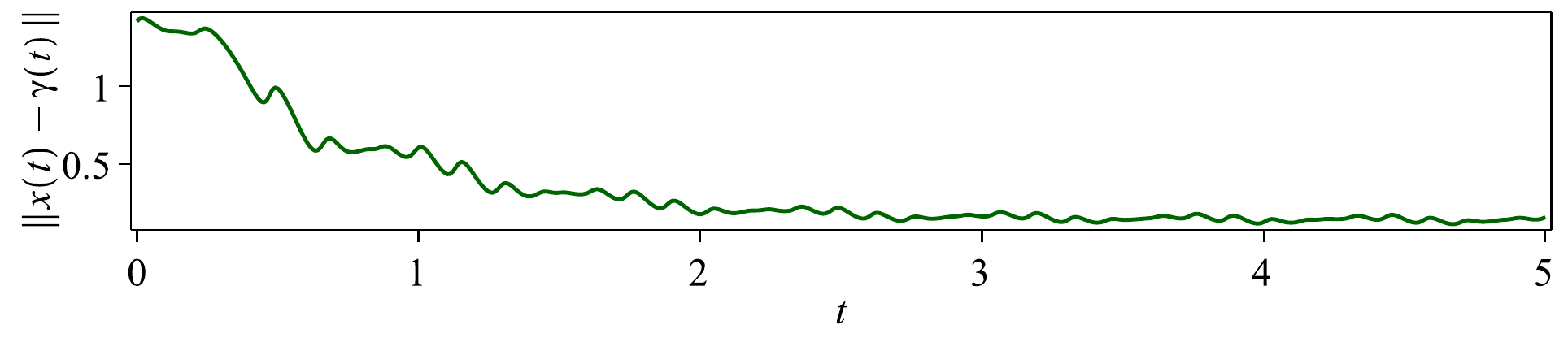}\\
      \caption{\footnotesize Trajectory (top), controls and the tracking error $\|x(t)-\gamma(t)\|$ (bottom)  for system~\eqref{uni} with controls~\eqref{cont1} (left),~\eqref{cont2} (center) and~\eqref{cont4} (right).~\label{figJ}}}
   \end{figure*}

Let $\Omega=5$, $\omega=50$, and
$
\gamma(t)=(0.1t,\sin(0.1t))^T,
$
so that the cost function is of the form
\begin{equation}\label{J1}
 \hat J(x,\gamma)= J(x-\gamma)=(x_1-0.1t)^2+{(x_2-\sin(0.1t))^2}.
\end{equation}
{In all simulations, we assume that system \eqref{uni} is initialized at $ x_1(0) = -1 ,x_2(0) = 1 $,  {the functions $F_1,F_2$ satisfy~\eqref{class}, and $\vartheta=1$}.}
{For the first case, take} 
\begin{equation}\label{cont1}
u(t)=\sqrt{\alpha\omega}\big(J(x-\gamma){\cos(\omega t)}+{\sin(\omega t)} \big).
\end{equation}
Here and in the sequel,  we denote $u(t):=u^\omega(t,J(x-\gamma))$.
Such type of controls were introduced in~\cite{DurrAuto} {{and also used in other classical extremum seeking approaches (e.g., \cite{Kr00})}, possibly with additional filters}.
{{The main advantages of this control are its simple analytical form and applicability for a wide class of cost functions.}}
The corresponding plots are shown in  Fig.~\ref{figJ} (left).
The {following} control, { introduced  in~\cite{Sch14}, possesses similar properties (and, moreover, has an a priori known bound):}
\begin{align}
&u(t)=\sqrt{\frac{\alpha\omega}{2}}\big(\sin(J(x-\gamma)){\cos(\omega t)}+\cos(J(x-\gamma)){\sin(\omega t)} \big),\nonumber\\
&|u(t)|\le\sqrt{{\alpha\omega}} \quad \text{ for any }J \; \text{and all }x\in\mathbb R^2,\,t\ge 0. \label{cont2}
\end{align}
The corresponding plots are shown in  Fig.~\ref{figJ} (center).

Although the controls~\eqref{cont1} and~\eqref{cont2} possess several useful properties, they always
lead to non-vanishing  oscillations of the trajectories of the extremum seeking system. This can be explained, in particular, by the fact that the controls~\eqref{cont1} and~\eqref{cont2} do not vanish for $J=0$, i.e.,
when approaching the target. Requiring $F_1(J),F_2(J)\to 0$ as $J\to 0$, it is possible to construct control laws which reduce the amplitude of oscillations and ensure better convergence properties.  In particular, these properties are ensured with the following control {proposed in~\cite{SD17}}:
\begin{equation}\label{cont3}
\begin{aligned}
u(t)=&\sqrt{\alpha\omega}\big(\sqrt{ |J(x-\gamma)|}\sin(\ln(|J(x-\gamma)|)){\cos(\omega t)}\\
&+\sqrt{|J(x-\gamma)|}\cos(\ln(|J(x-\gamma)|)){\sin(\omega t)}\big).
\end{aligned}
\end{equation}
%
In order to combine the advantages of controls~\eqref{cont3} (i.e., vanishing amplitudes when reaching the target) {and~\eqref{cont2}} (i.e., bounded excitation signals independent of the cost function),
the following control function has been proposed in~\cite{GZE}:
\begin{equation}\label{cont4}
u(t)=\sqrt{\alpha\omega\phi}\big(\sin(\psi){\cos(\omega t)}+\cos(\psi)){\sin(\omega t)} \big),
\end{equation}
  where $\phi={\frac{1{-}e^{{-}|J(x-\gamma)|}}{1{+}e^{J(x-\gamma )}}}$, $\psi=e^{|J(x-\gamma)|}{+}2\ln(e^{|J(x-\gamma)|}{-}1)$ for $J{\ne} 0$, and $u{=}0$ for $J{=}0$.
  Fig.~\ref{figJ} (right) presents the simulation results for system~\eqref{uni} with the controls~\eqref{cont4}. It can be seen that the controls~\eqref{cont3} and~\eqref{cont4} exhibit smaller tracking error and the control amplitudes.
{However, it has to be noticed that both controls~\eqref{cont3} and~\eqref{cont4} exhibit better behavior of an extremum seeking system only in case of known minimal value of the cost function, and control~\eqref{cont4} requires that $x(t_0)$ is close enough to $x^*(\gamma(t_0))$ (under a proper scaling of the cost function) for better convergence properties.}

\subsection{Experimental results}

\begin{figure}[b]
	\centering
	\includegraphics[width=0.4\linewidth]{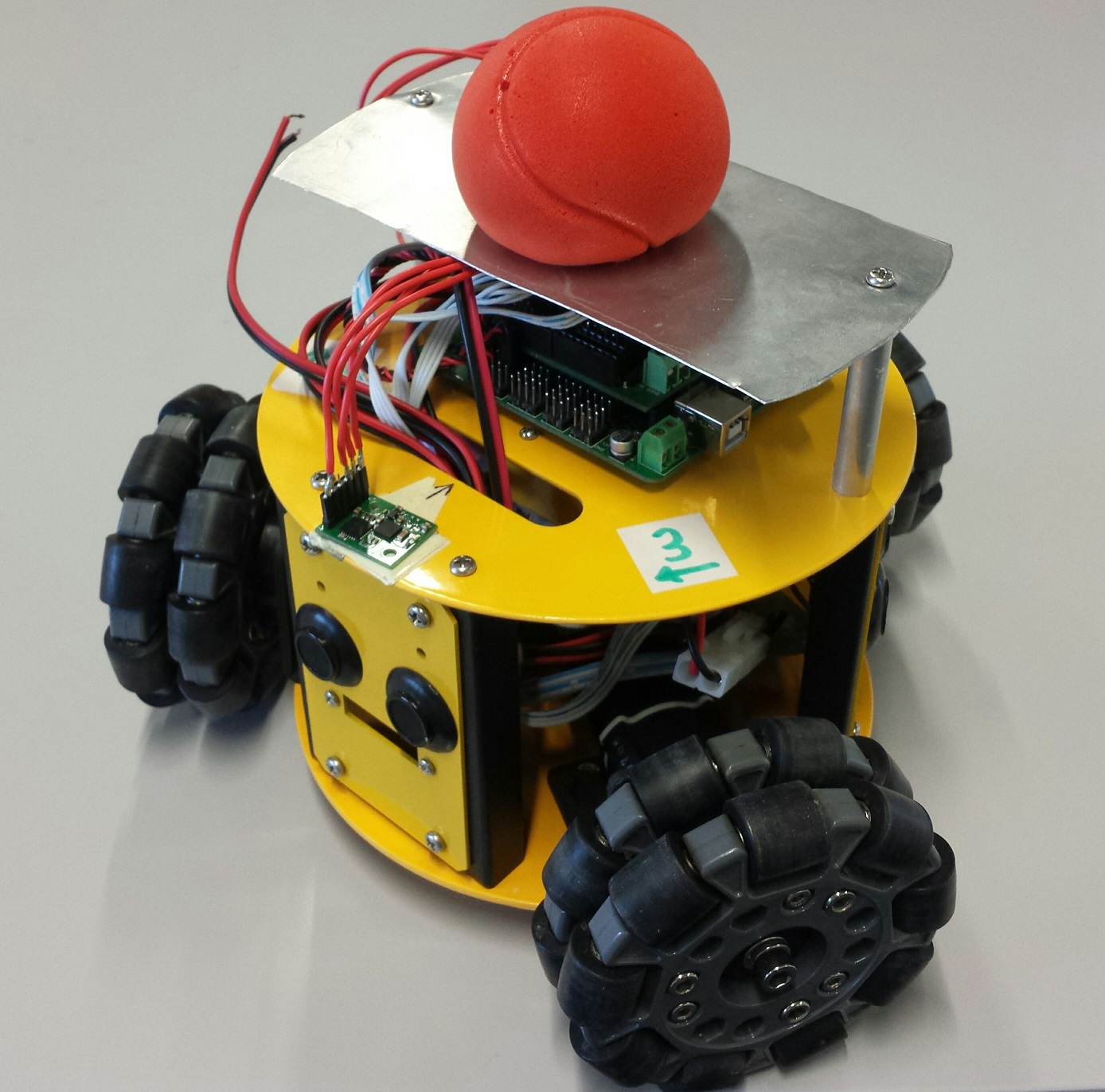}
	\caption{\footnotesize The omni wheel robot used in the experiments.}\label{figRobot}
\end{figure}

The above examples show that the proposed new extremum seeking
control laws perform very well in numerical simulations.
In this section we want to illustrate that the benefits also
transfer to an experimental setup. We validated the control
on a three-wheeled mobile robot (see Fig.~\ref{figRobot})
both in a fixed and a moving target tracking scenario.

Due to limitations in the experimental setup we do not directly
measure the distance to the target but evaluate it using $ (x_1,x_2) $-position measurements of both the robot and the target obtained
from tracking them with a camera.

In the fixed target scenario, we let $ \omega = 3 $, $ \Omega = 1.5 $
and assume the cost function to be of the form~\eqref{J}
with $ \gamma \equiv x^\star = ( 0.5, 0.7 )^\top $ being the constant position of
the target. 
We compared the control laws~\eqref{cont1}, \eqref{cont2}
and \eqref{cont4} where the parameters $ \alpha $ and $ \kappa $ were tuned
under the assumption that the input is bounded as $ \vert u(t) \vert \leq 0.4 $,
see Table~\ref{tabParams}.
\begin{table}[h]
	\centering
	\begin{tabular}{cllll}
		\toprule
		\textbf{Control}       & \multicolumn{2}{c}{\textbf{Fixed target}} & \multicolumn{2}{c}{\textbf{Moving target}} \\
		\textbf{law}           & $ \alpha $              & $ \kappa $ & $ \alpha $            & $ \kappa $ \\ \midrule
		\eqref{cont1} & $2.25\cdot 10^{-4}$     & $10$       & --                     & -- \\[0.5em]
		\eqref{cont2} & $4.84\cdot 10^{-2}$     & $4$        & $ 5.29\cdot 10^{-2} $ & $4$ \\[0.5em]
		\eqref{cont4} & $3.249\cdot 10^{-1}$    & $4$        & $ 2.5 \cdot 10^{-1} $ & $1$\\
		\bottomrule
	\end{tabular}
	\caption{\footnotesize Parameters used in the experiments.}\label{tabParams}
\end{table}

\begin{figure}
	\centering
	\includegraphics[width=1\linewidth]{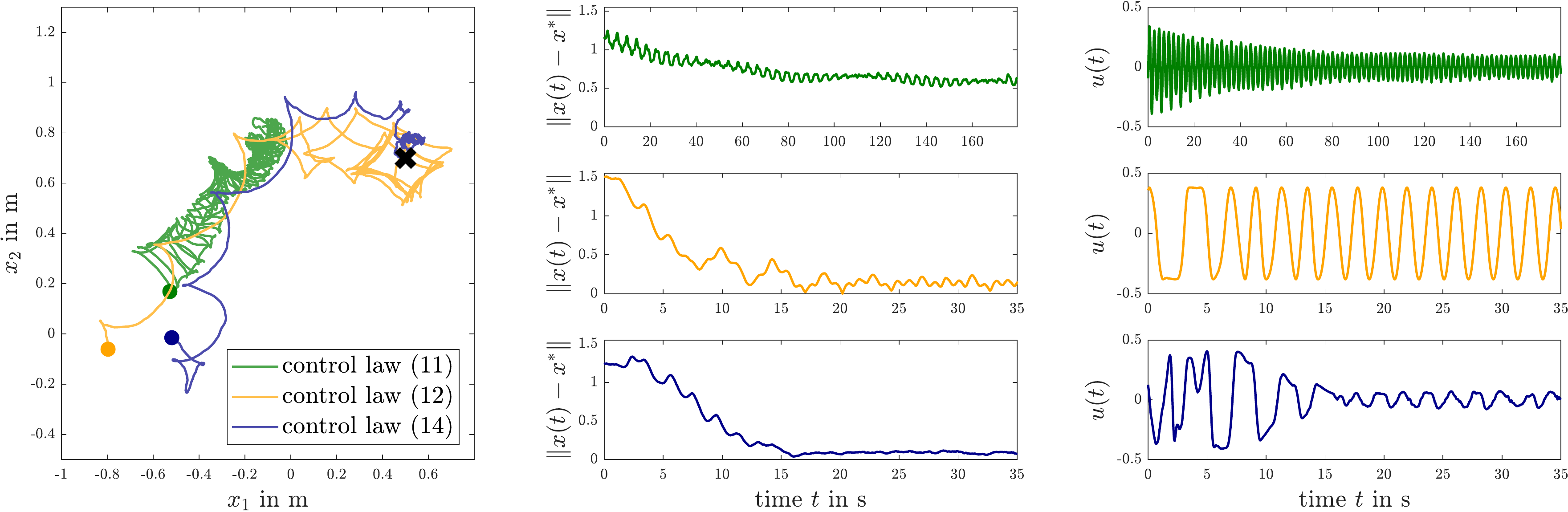}
	\caption{\footnotesize Experimental results for the fixed target scenario.
	The plot on the left-hand side shows the trajectories of the
	robot in the $(x_1,x_2)$-plane using the three different
	control laws~\eqref{cont1} (green), \eqref{cont2} (orange) and \eqref{cont4} (blue).
	The target point is indicated by a black cross. The plots
	in the middle depict the evolution of the distance
	to the target, again for the control laws~\eqref{cont1} (green, upper), \eqref{cont2} (orange, middle) and \eqref{cont4} (blue, lower).
	The plots on the right-hand side depict the control input
	$ u(t) $.
	\label{figFixedTarget}}
	\centering
	\includegraphics[width=0.49\linewidth]{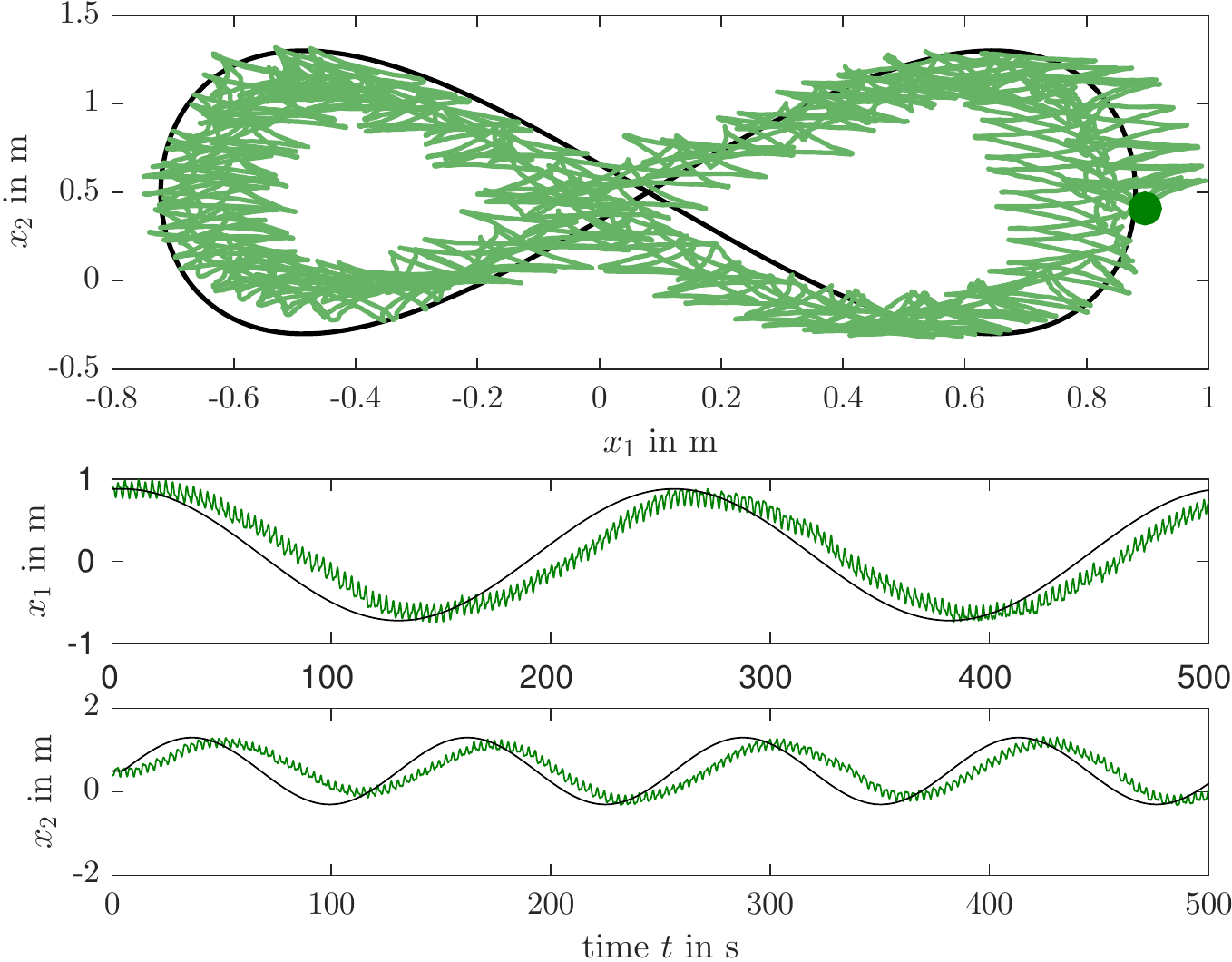}
	\includegraphics[width=0.49\linewidth]{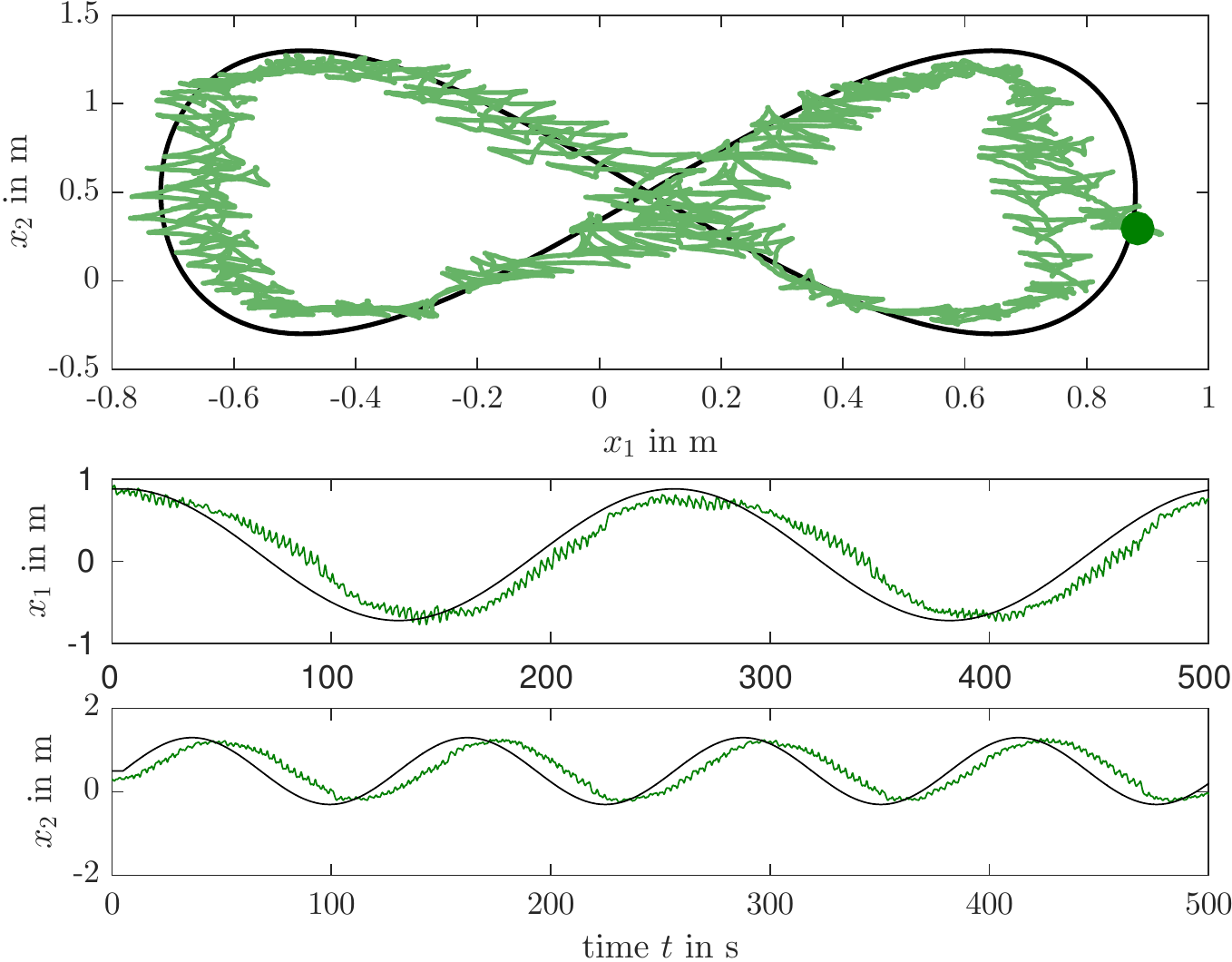}
	\caption{{\footnotesize Experimental results for the moving target scenario using control
	law \eqref{cont2} (left) and \eqref{cont4} (right). The target trajectory is depicted in black and the robot
	trajectory is depicted in green. The accumulated tracking error is
	$ \int_{0}^{500} \Vert x(t) - \gamma(t) \Vert^2 \mathrm{d}t \approx 367.7589 $ for \eqref{cont2} and $\approx 75.3409$ for \eqref{cont4},
	and the accumulated control {effort} is $ \int_{0}^{500} \vert u(t) \vert \mathrm{d}t \approx 129.2032 $ for \eqref{cont2} and $\approx 64.7733$  for \eqref{cont4}.}}\label{figMovingBounded}
\end{figure}

The experimental results are depicted in Fig. \ref{figFixedTarget}.
Control law~\eqref{cont1} shows the worst performance and does not
converge very close to the target, even in much longer time.
The performance of \eqref{cont2} and \eqref{cont4} is comparable
in terms of the accumulated squared distance error and the convergence
time. However, while control law \eqref{cont2} is non-vanishing and
thus the robot circulates around the target in the end, the robot
only makes small movements in the end when using control law \eqref{cont4},
and the resulting total control {effort} is drastically reduced in
comparison to control law \eqref{cont2}. The reason why the control
input does not vanish completely is the imperfect rotational motion
of the robot when $ u(t) = 0 $.

In the moving target scenario, we let $ \omega = 3 $, $ \Omega = 1 $.
The goal is to track a target moving along a figure eight curve,
i.e., the cost function takes the form~\eqref{J} with
\begin{align}
	\gamma = \big( 0.8 \cos( 0.025 t ) + 0.08, 0.8 \sin(0.05 t) + 0.5 \big)^\top .
\end{align}

We compared the control laws \eqref{cont2} and \eqref{cont4}.
Again, the parameters $ \alpha $ and $ \kappa $ were tuned
under the assumption that the input is bounded as $ \vert u(t) \vert \leq 0.4 $,
see Table~\ref{tabParams}. The experimental results are depicted in
Fig. \ref{figMovingBounded} (left) for control~\eqref{cont2} and in Fig. \ref{figMovingBounded} (right)
for control law \eqref{cont4}. Both control laws achieve tracking the moving target, where control law \eqref{cont4} shows a better
behavior in terms of the tracking error
while requiring only approximately  half the control {effort}.

All in all, the experimental results show that the new
extremum seeking control laws can lead to improved performance
also in practical implementations. Nevertheless,
due to low upper limits for $ \omega $ and $ \Omega $,
there is still quite a gap between experimental and simulative
results.

\section*{Conclusions}

In this paper, a novel family of extremum seeking laws have been introduced for unicycle models.
We have proved that the proposed controls can be utilized for tracking an
extremum point of a time-varying cost function by extending the theoretical results from \cite{GZE} and \cite{GDEZ17}.
{In particular,  we have discussed how the results
can be applied to moving target tracking problems.}
{We have illustrated by simulations as well as experiments} that the proposed family of extremum seeking laws
performs remarkably well for these type of problems.
Our next goals are to construct further extensions of the family of control functions~\eqref{class} for more general classes of cost functions and to evaluate their performance with simulations and
experiments.

\bibliographystyle{plain}
\bibliography{biblio_es}

\section*{{Appendix A. Proof of Theorem~\ref{thm_quad}}}

\textit{Step 0. Preliminary constructions.}\\
Let $\lambda\in(0,\rho)$, $\vartheta>\bar \vartheta=\frac{\nu}{2\sqrt{\kappa\lambda}}$, $\delta\in\Big(0,\frac{\sqrt\rho-\sqrt\lambda}{\sqrt\kappa}\Big)$, $x^0\in B_{\delta}(\mathcal L_{\lambda,t_0})$.
Introducing the new variables
  \begin{equation}\label{zam}
  \xi= x-x^*(\gamma),
  \end{equation}
we rewrite system~\eqref{uni} with controls~\eqref{cont}  as 
\begin{equation}\label{afft}
  \dot \xi=-\dot\gamma+\sum_{j=1}^4g_j(J(\xi))\sqrt\omega v_j(\omega t),
\end{equation}
with $J(\xi)=\kappa\|\xi\|^2$, $
g_i(J(\xi)){=}\big(F_i(J(\xi)),0\big)^T,$
$g_{i+2}(J(\xi)){ =}\big(0,F_i(J(\xi))\big)^T,\ i=1,2$, $v_i$ defined as in the proof of Theorem~2, and $\xi(t_0)=\xi^0=x(t_0)-\gamma(t_0)$,
$$\xi^0\in \widetilde D_0=B_{\delta}\big(\{\xi\in\mathbb R^2:\kappa\|\xi\|^2\le\lambda\}\big)=B_{\delta+\sqrt{{\lambda}/{\kappa}}}(0).$$
We denote  $\widetilde D=\{\xi\in\mathbb R^2:J(\xi)\le\rho\}$, take $K>1$ such that $\vartheta>K\bar \vartheta$,   and fix any $\mu,\delta_0,\delta_{\min}$ satisfying
$$1<\mu<K, \,\frac{\mu^2\lambda}{K^2}<\delta_0<\lambda, \,\delta_{\min}\in(0,\delta_0).$$
Since $F_i\circ J\in C(\widetilde D;\mathbb R^2)$,
there is a $\tau_0>0$ such  that the solutions of~\eqref{afft} are well-defined in $\widetilde D$ for all $t\in[t_0;t_0+\tau_0]$.\\
\textit{Step 1. A priori bounds of the solutions.}
Consider the function $w(t)=\|\xi(t)-\xi^0\|$. Estimating the derivative  of $w^2(t)$ along the trajectories of~\eqref{afft} with regard to the assumptions of this theorem, we get
$$
\begin{aligned}
\dot w(t)&\le \|\dot\gamma(t)\|+\sqrt\omega\sum_{j=1}^4\|g_j(J(\xi))\|| v_j(\omega t)|\\
&\le  \nu+ M\sqrt\omega \quad \text{ with } M=2\sqrt{\alpha \vartheta}\max\limits_{\xi\in \widetilde D}{\sum\limits_{i=1,2}}\Big|F_i(J(\xi))|.
\end{aligned}
$$
Solving the corresponding comparison equation with $w(t_0)=0$, we conclude that $w(t)\le (\nu+ M\sqrt\omega)(t-t_0)$.
Hence, for all $t\in[t_0;t_0+T]$ ($T=\tilde k/\omega$, $\tilde k=2\pi k$), we have $\|\xi(t)-\xi^0\|\le \frac{\tilde k(\nu+M\sqrt{\omega})}{\omega}$.
Define $$d=\min\Big\{\frac{\sqrt\lambda-\sqrt\delta_0}{\sqrt\kappa},
\frac{\sqrt\delta_0-\sqrt\delta_{\min}}{\kappa},\frac{\sqrt\rho-\sqrt\lambda}{\sqrt\kappa}-\delta\Big\}.$$
Then,
for all  $\omega>\omega_1=\max\Big\{\frac{\tilde k}{\tau_0},\frac{4\nu\tilde k^2}{(\sqrt{M\tilde k^2+4d\nu\tilde k}-M\tilde k)^2}\Big\}$, $t\in[t_0;t_0+T]$, the following properties hold:\\
(P1) $\|\xi(t){-}\xi^0\|< d<\frac{\sqrt\lambda}{\sqrt\kappa}$;\\
(P2) $\xi^0\in \widetilde D_0\Rightarrow \xi(t)\in\widetilde D$ for all   $t\in[t_0;t_0+T]$;\\
(P3) $J(\xi^0)\le\delta_0\Rightarrow J(\xi(t))<\lambda$ for all   $t\in[t_0;t_0+T]$;\\
(P4) $J(\xi^0)>\delta_0\Rightarrow J(\xi(t))>0$ for all   $t\in[t_0;t_0+T]$.\\
\textit{Step 2. Representation of the solutions.} Let us expand the solutions of system~\eqref{afft} into the Volterra-type series. From (P4) and (B1),
the  representation
\begin{align*}
g_j(J(\xi(t)))&=g_j(J(\xi^0))+\int_{t_0}^t\frac{dg_j(J(\xi(s)))}{ds}ds\\
&=g_j(J(\xi^0))+\int_{t_0}^t\frac{\partial g_j(J(\xi(s)))}{\partial \xi}\dot\xi(s)ds=g_j(J(\xi^0))\\
&{-}\int\limits_0^t\Big(L_{\dot\gamma}(g_j\circ J)(\xi(s))+\sum_{k=1}^4L_{g_k}(g_j\circ J)(\xi(s))\sqrt\omega v_k(\omega s)\Big)ds
\end{align*}
is well-defined for all   $t\in[t_0;t_0+T]$.
 Applying the same procedure to $L_{F_j}(F_i\circ J)(\xi(s))$ and using
 $$\xi(t){=}\xi^0{-}{\int_0^t}\dot\gamma(\xi(\tau))d\tau{+}{\sum\limits_{j=1}^4}{\int_0^t}g_j(J(\xi(\tau)))\sqrt\omega v_j(\omega\tau)d\tau,$$
 we get
\begin{align}
\xi(t)&=\xi^0-\int\limits_{t_0}^t \dot\gamma(\tau)d\tau{+}\sum_{j=1}^4g_j(J(\xi^0))\sqrt\omega\int\limits_{t_0}^t v_j(\omega \tau)d\tau \label{volt}\\
&{+}\sum_{j,k=1}^4L_{g_k}(g_j\circ J)(\xi^0)\omega\int\limits_{t_0}^t\int\limits_{t_0}^\tau v_k(\omega s)v_j(\omega\tau)dsd\tau+r(t),\nonumber
\end{align}
where $r(t)$ is the remainder,
\begin{align*}
r(t)&=\omega^{3/2}\sum_{j,k,l=1}^4\int\limits_{t_0}^t\int\limits_{t_0}^\tau\int\limits_{t_0}^s L_{g_l}L_{g_k}(g_j\circ J)(\xi(p))v_l(\omega p) v_k(\omega s)\nonumber\\
  &\qquad\qquad\times v_j(\omega\tau)dpdsd\tau\\
 &+\sqrt\omega \sum_{j=1}^4\int\limits_{t_0}^t\int\limits_{t_0}^\tau L_{\dot\gamma}(g_j\circ J)(\xi(s))v_j(\omega\tau)dsd\tau\nonumber\\
  &+\omega\sum_{j,k=1}^4\int\limits_{t_0}^t\int\limits_{t_0}^\tau\int\limits_{t_0}^sL_{\dot\gamma}L_{g_k}(g_j\circ J)(\xi(p)) v_j(\omega\tau)v_k(\omega s)dpdsd\tau.\nonumber
\end{align*}
 It can be shown that
$$\|r(t_0+T)\|\le \tilde k^2\omega^{-3/2}R,$$
where
$$
\begin{aligned}
R=&\frac{\tilde k(H+\nu M_2\omega^{-1/2})}{3}+\nu M_1, \,M_1{=}\max\limits_{\xi\in \widetilde D}{\sum\limits_{i=1,2}}\Big\|\nabla F_i(J(\xi))\Big\|,\\
&M_2{=}\max\limits_{\xi\in \widetilde D}{\sum\limits_{i,j=1,2}}\Big\|\nabla\big(L_{F_j}(F_i\circ J)\big)(\xi)\Big\|,\\
&H=\max\limits_{\xi\in \widetilde D}\sum\limits_{i=1,2}\Big\|L_{F_k}L_{F_j}(F_i\circ J)(\xi)\Big\|.
\end{aligned}
$$
Recall that the above represenation of the solutuons of~\eqref{afft} is well-defined for all   $t\in[t_0;t_0+T]$ because of (P4) and (B1).\\
\textit{Step 3. Estimation of the cost function.} Direct calculation of integrals in~\eqref{volt} for $t=T+t_0$ and the application of formula~\eqref{class} imply
\begin{align*}
\xi(t_0+T)&=\xi^0-\int\limits_{t_0}^t \dot\gamma(\tau)d\tau-T (\vartheta\nabla J(\xi^0)^T-\Phi(J(\xi^0)))+r(t)\\
&\le \xi^0-T \big(\vartheta\nabla J(\xi^0)^T-\nu-\Phi(J(\xi^0))\big)+\tilde k^2\omega^{-3/2}R,
\end{align*}
where  $\Phi(J(\xi))=(-\varphi_2(\xi),\varphi_1(\xi))^T$,
 \begin{align*}
 \varphi_s(\xi)=& \frac{1}{2k}\nabla_{\xi_s}\Big(F_1^2( J(\xi))+F_2^2( J(\xi)\Big),\;s=1,2.
\end{align*}
Note that $\big(\Phi(J(\xi)),\nabla J(\xi)\big)=0$ for each $\xi$.\\
 Define $y=\xi(t_0+T)-\xi^0$, $\eta{=}(1{-}\theta)\xi^0{+}\theta \xi(t_0+T)$, $\theta{\in}(0,1)$.
Using the Taylor expansion of $J$ with the Lagrange form of the remainder,
  \begin{align*}
       J(\xi(t_0+T))=J(\xi^0)+{\nabla J(\xi^0)y}+\frac{1}{2}\sum_{i,j=1}^n\frac{\partial^2 J(\xi)}{\partial \xi_i\partial \xi_j}\bigg|_{\xi=\eta}y_iy_j,
   \end{align*}
we get the following estimate for  $J(x^0)> \delta_0$:
    \begin{align*}
J(\xi(t_0+T))&\le J(\xi^0)-T \vartheta\|\nabla J(\xi^0)\|^2\\
       &+T\big(\Phi(J(\xi^0)),\nabla J(\xi^0)\big)+\|\nabla J(\xi^0)\|(T\nu+\tilde k^2\omega^{-3/2}R)\\
       &+\kappa\tilde k^2\omega^{-2}\Big(\vartheta\|\nabla J(\xi^0)\|+\nu+\|\Phi(J(\xi^0))\|+\omega^{-1/2}R\Big)^2\\
       &\le J(\xi^0)-\|\nabla J(\xi^0)\|^2T\big( \vartheta-\frac{\nu}{2\sqrt{\kappa\delta_0}}-\tilde k^2\omega^{-3/2}\widetilde R\big),\\
      &\widetilde R=\frac{\nu R(T)}{2\sqrt{\kappa\delta_0}}+\kappa\omega_1^{-1/2}\Big(\vartheta+L+\frac{\omega_1^{-1/2} R}{2\sqrt{\kappa\delta_0}}\Big)^2>0,
           \end{align*}
 where $L=\max\limits_{\xi\in \widetilde D}\sum\limits_{i=1,2}\Big\|L_{F_i}(F_i\circ J)(\xi)\Big\|$.\\
By the definition of $\delta_0$, $\vartheta-\frac{\nu}{2\sqrt{\kappa\delta_0}}>(1-{\mu^{-1}})\bar \vartheta=:\bar\beta>0.$
For any $\beta\in(0,\bar \beta)$, let $\omega_2=\max\Big\{\omega_1,\Big(\frac{\tilde k^2 \widetilde R}{\bar\beta-\beta}\Big)^{2/3}\Big\}$. \\
Then, for all $\omega>\omega_2$, $\xi_0\in\widetilde D_0\setminus\{\xi:J(\xi)\le\delta_0\}$,
$$
\begin{aligned}
J(\xi(t_0+T))&\le J(\xi^0)-4\kappa TJ(\xi^0) \big(\bar\beta-\tilde k^2\omega^{-3/2}\widetilde R\big)\\
&\le  J(\xi^0)-4\kappa T J(\xi^0) \beta.
\end{aligned}
$$
Defining $\omega_3=\max\{\omega_2,4\kappa\tilde k \beta\}$, we conclude that $4\kappa \tilde k\omega^{-1} \beta<1$ for all $\omega>\omega_3$. Therefore, $\xi(T)\in\widetilde D_0$, and  the last estimate can be rewritten as
\begin{equation}\label{estJ}
J(\xi(t_0+T))\le J(\xi^0)e^{-4\kappa \beta T}.
\end{equation}
\textit{Step 4. Attractivity.} On this step we show that there exists an  $N\in\mathbb N\cup\{0\}$ such that $J(\xi(NT))<\delta_0$, and $J(\xi(t))\le\lambda$ for all $t\ge NT$.
Suppose that $J(\xi(pT))\ge\delta_0$, for all $p=0,1,2,\dots$, and take  $N=\Big[\frac{1}{4\kappa \beta T}\ln\frac{(\delta+\sqrt{\lambda/\kappa})^2}{\delta_0}\Big]+1$. Then the iteration of~\eqref{estJ} with $t=t_0+T,t_0+2T,\dots$ gives
 $$
J(\xi(t_0+N T))\le J(\xi^0)e^{-4\kappa \beta N T}<\delta_0.
$$
So, we get the contradiction which proves that there exists an $N>0$ such that  $J(\xi(NT))<\delta_0$. Thus, we have two possibilities. If  $J(\xi(t))<\delta_0<\lambda$ for all $t\ge NT$, then the proof of the attractivity is completed.
Otherwise, we recall from (P3) that $J(\xi(t))<\lambda$ for all $t\in[NT,(N+1)T]$.
This again  yields two possibilities:\\
a)  $J(\xi((N+1)T))<\delta_0$;\\
b) $\delta_0\le J(\xi((N+1)T))<\lambda$, so that we can apply estimate~\eqref{estJ}.
Repeating the above argumentation, we obtain $J(\xi(t))\le\lambda$ for all $t\ge NT$.\\
\textit{Step 5. Decay rate.}
Without loss of generality, assume that $J(\xi(pT))\ge\delta_0$ for all $p=\overline{0,N-1}$. Then
\begin{equation}\label{estJ2}
J(\xi(t))\le J(\xi^0)e^{-4\kappa \beta (t-t_0)}\quad \text{ for }t=t_0+T,\dots, t_0+NT.
\end{equation}
The estimate~\eqref{estJ2} together with (P2), (P3), and the results of Step~4 implies that the solutions of system~\eqref{afft} with the initial conditions from $\widetilde D_0$ are well-defined in $\widetilde D$ for all $t\ge t_0$.
It remains to estimate $\|\xi(t)\|$ for the solutions of~\eqref{afft} if $t\in [t_0,t_0+NT]$.
For any $t\in[t_0,t_0+NT]$, we
denote the integer part of $tT$ as $t^T_{int}$, and observe that $0<t-{t^T_{int}}{T}<T$.
Using~(P1), we obtain that, for all $t\in[t_0,t_0+NT]$,
\begin{align*}
J(\xi(t))
&\le\Big(J^{1/2}(\xi(t^T_{int}T))+\sqrt\kappa\|\xi(t)-\xi(t^T_{int}T)\|\Big)^2\le\Big(J^{1/2}(\xi^0)e^{-2\kappa  \beta (t^T_{int}T-t_0)}+\sqrt\lambda\Big)^2\\
&\le \Big(e^{2\kappa  \beta T}J^{1/2}(\xi^0)e^{-2\kappa  \beta (t-t_0)}+\sqrt\lambda\Big)^2.
\end{align*}
Formula~\eqref{zam} completes the proof:
 for  $\lambda\in(0,\rho)$, $\delta\in\Big(0,\frac{\sqrt\rho-\sqrt\lambda}{\sqrt\kappa}\Big)$, $\vartheta>\frac{\nu}{2\sqrt{\kappa\lambda}}$, we may take $\omega_0>\omega_3$, $\tau(\omega_0)>\Big[\frac{1}{4\kappa \omega_0 \beta \tilde k}\ln\frac{(\delta+\sqrt{\lambda/\kappa})^2}{\lambda}\Big]+1$, and conclude that, for all  $x^0\in B_{\delta}(\mathcal L_{\lambda,t_0})$, $\omega>\omega_0$,  the solutions of system~\eqref{uni} with controls~\eqref{cont} satisfy the following property:
\begin{align*}
&\|x(t)-\gamma(t)\|\le \Big(e^{\frac{4\pi k\kappa  \beta}{\omega}}\|x^0-\gamma(t_0)\|e^{-2\kappa  \beta (t-t_0)}+\sqrt{\frac{\lambda}{\kappa}}\Big),\\
&\qquad\qquad\qquad\qquad\qquad\qquad\qquad\text{ for }t< t_0+\tau,\\
&x(t)\in\mathcal L_{\lambda,t}\;\text{ for }t\in[ t_0+\tau,\infty).
\end{align*}


\end{document}